\documentclass[12pt,twoside,reqno]{amsart}

\usepackage{amsthm, amsmath, amscd, amssymb,centernot}
\usepackage[left=2.5cm,top=2.5cm,bottom=3cm,right=2.5cm]{geometry}
\usepackage[pagebackref]{hyperref}

\newtheorem*{rep@theorem}{\rep@title}
\newcommand{\newreptheorem}[2]{%
	\newenvironment{rep#1}[1]{%
		\def\rep@title{#2 \ref{##1}} 
		\begin{rep@theorem}}%
		{\end{rep@theorem}}}
	\makeatother
	
	\newreptheorem{theorem}{Theorem}

\newcommand{\proj}{\mathbb{P}}

\newcommand{\Xr}{\widetilde{X}_r(x)}

\theoremstyle{plain}
\numberwithin{equation}{section}
\newtheorem{theorem}{Theorem}[section]
\newtheorem*{theorem*}{Theorem}
\newtheorem{proposition}[theorem]{Proposition}
\newtheorem{lemma}[theorem]{Lemma}
\newtheorem{corollary}[theorem]{Corollary}
\newtheorem{conjecture}[theorem]{Conjecture}
\newtheorem{definition}[theorem]{Definition}
\theoremstyle{definition}
\newtheorem{question}[theorem]{Question}
\newtheorem{remark}[theorem]{Remark}
\newtheorem{example}[theorem]{Example}

\begin{document}
	\title{Lower bounds for Seshadri constants on blow Ups of $\proj ^2$ }
	
	\author[Cyril J. Jacob]{Cyril J. Jacob}
	\address{Chennai Mathematical Institute, H1 SIPCOT IT Park, Siruseri, Kelambakkam 603103, India}
	\email{cyril@cmi.ac.in}
	
	\begin{abstract}
		Let $\pi: X_r \rightarrow \proj^2$ be a blow up of $\proj^2$ at $r$ distinct points $p_1,p_2,\dots, p_r$. We study lower bounds for Seshadri constants of ample line bundles on $X_r$. First, we consider the case when the points lie on a curve of degree $d\le 3$, and the case when $r\le 8$. We then assume that the points are very general and show that $\varepsilon(X_r)\geq \frac{1}{2}$ if the Strong SHGH conjecture is true. 
	\end{abstract}
	\subjclass[2020]{14C20, 14E05, 14H50}
	\keywords{Seshadri constants}

	\date{\today}

	\maketitle
	\section{Introduction}
	Throughout this article, we work over the field of complex numbers.
	
	Seshadri constants are numerical invariants associated with an ample line bundle on a projective variety. They measure the local positivity of the line bundle and were introduced by Demailly in \cite{Dema} to study the Fujita conjecture. 
	\begin{definition}
		Let $X$ be a smooth complex projective surface. The \textit{Seshadri constant} of a nef line bundle $L$ on a surface $X$ at a point $x \in X$ is defined as:
		$$\varepsilon(X,L,x)=\inf\bigg \{ \frac{L\cdot  C}{\text{mult}_x C}~~\bigg  |  ~~C \text{ is a reduced irreducible curve on $X$ passing through $x$} \bigg  \}. $$
	\end{definition}

	We also define three related constants on a surface $X$, each connected to the Seshadri constant of an ample line bundle at a point. First, the \textit{Seshadri constant of $X$ at a point $x\in X$}, denoted by $\varepsilon(X,x)$, is defined as follows:
	$$\varepsilon(X,x)=\inf_{L\text{ ample }}\varepsilon(X,L,x).$$
	
	Second, for a fixed ample line bundle $L$, the \textit{Seshadri constant of  $X$ at $L$ } is defined as:
	$$\varepsilon(X,L)=\inf_{x\in X} \varepsilon(X,L,x).$$
	
	Finally, we define the \textit{Seshadri constant of $X$}, $\varepsilon(X)$ as:
	$$\varepsilon(X)=\inf_{x\in X} \varepsilon(X,x)=\inf_{L \text { ample}} \varepsilon(X,L).$$
	
 More details about Seshadri constants and local positivity can be found in \cite[Section 5]{Laz} and \cite{Bau}.
	
	Computing exact values of Seshadri constants is difficult in general, and remains open even for blow-ups of \( \mathbb{P}^2 \). As a result, a central problem in this area is to determine effective lower bounds. A seminal result in this direction is due to Ein and Lazarsfeld \cite{EinLaz}, who proved that \( \varepsilon(X, x) \geq 1 \) for very general points \( x \in X \). On the other hand, Miranda constructed examples where the Seshadri constant can be arbitrarily small (see \cite[Example 3.1]{EinLaz} and \cite[Example 5.2.1]{Laz}). These contrasting results motivate the following natural and fundamental question:
	
	\begin{question}
		Does there exist a surface $X$ such that $\varepsilon(X)=0$?
	\end{question} 
	In this article, we will address this question when $X$ is a blow up of $\proj^2$ at finitely many points. Let $\pi:X_r\to \mathbb{P}^2$ be the blow up of $\mathbb{P}^2$ at distinct points $p_1,p_2,\dots,p_r$. We focus on obtaining lower bounds for Seshadri constants on $X_r$, under various configurations of $p_1,p_2,\dots,p_r$.
	
	We obtain a lower bound for $X_r$ when $p_1,p_2,\dots,p_r$ lie on a curve of degree $d\le 3$ as follows:
		\begin{reptheorem}{dcurve}
		Let $p_1,p_2,\dots,p_r$ be r distinct points lying on a curve $C\subseteq \mathbb{P}^2$ of degree $d\leq 3$. Let $X_r$  denote the blow up of $\mathbb{P}^2$ at $p_1,p_2,\dots ,p_r$. Then for all line bundle bundles $L$ and for all points $x\in X_r$ we have,
		$$
		\varepsilon(X_r,L,x) \geq 
		\begin{cases}
			1 \text{ if $d=1,2$ or $d=3$ and $L\cdot (-K_{X_r})\geq 2$}\\
			\frac{1}{2} \text{ if $d=3$ and $L\cdot (-K_{X_r})=1.$}
		\end{cases}
		$$
	\end{reptheorem}
	
	Then we study the case $r\le 8$ and show that $\varepsilon(X_r)\ge \frac{1}{2}$. More precisely, by Theorem \ref{r7geq1} and Theorem \ref{X8} we have the following:
	\begin{theorem*}
		For $r\le 8$, let $X_r$ denote the blow up of $r$ distinct points on $\mathbb{P}^2$ and $L$ be an ample line bundle on $X_r$. If $\varepsilon(X_r,L,x)<1$, then
		\begin{enumerate}
			\item $r=8$,
			\item $p_1,p_2,\dots ,p_8$ are smooth points on an irreducible singular cubic,
			\item $L=-K_{X_8}$, and
			\item $\varepsilon(X_8,L,x)=\frac{1}{2}$.
		\end{enumerate}
	\end{theorem*}
	
		 Now to address the same for $r\ge 9$, we assume  a conjecture concerning linear systems of plane curve called Generalized SHGH or Strong SHGH (see Conjecture \ref{SSHGH}), which is a stronger version of the famous SHGH conjecture and also we assume the blown up points are very general. As stated above Miranda showed that Seshadri constants can be arbitrarily small.  In fact, Miranda's construction shows that for every integer $m\ge 2$, there is a blow up $X\to \mathbb{P}^2$ at $d^2$ points on a degree $d$ plane curve (the choice of $d$ depends on $m$), a particular ample line bundle $L$ on $X$ and a particular point $x\in X$ such that $\varepsilon(X,L,x)=\frac{1}{m}$. In contrast we show that if $X \to \mathbb{P}^2$ is a blow up of $\mathbb{P}^2$ at \textit{very general} points, then Conjecture \ref{SSHGH} implies that $\varepsilon(X,L,x)\ge \frac{1}{2}$ for \textit{every} ample line bundle $L$ on $X$ and \textit{every} point $x\in X$:

	 
	 	\begin{reptheorem}{scgeq1}
	 	Assume that Conjecture \ref{SSHGH} is true. Let $\pi : X_r \to \proj ^2$ be the blow up of $\proj ^2$ at very general points $p_1,p_2,\dots ,p_r$. Then 
	 	$$\varepsilon(X_r)\geq \frac{1}{2}.$$
	 	Further, if $\varepsilon(X_r,L,x)<1$ for some ample line bundle $L$ and for some $x\in X$, then $\varepsilon(X_r,L,x)=\frac{1}{2}$.
	 \end{reptheorem}

		\subsection*{Acknowledgments}
	It is a pleasure to thank my advisor Krishna Hanumanthu for a meticulous reading of the manuscript and for many helpful discussions, suggestions, and comments. I thank the referee for a careful reading of the paper and numerous suggestions, which
	improved the paper. The author was partially supported by a grant from Infosys Foundation.
	
	\section{Preliminaries}
	This section recalls the key definitions and conjectures necessary for our results. We begin with an alternative characterization of the Seshadri constant that will be used extensively throughout the paper.
	\begin{lemma}
			Let $X$ be a smooth projective surface and $x\in X$ be a point on $X$. Let $\pi_x : \widetilde{X}(x) \rightarrow X$ be the blow up of $X$ at $x$. Let $E=\pi_x^{-1}(\{x\})$ denote the exceptional divisor. Then
		$$\varepsilon(X,L,x)=\sup\left\{t\in \mathbb{R}_{\ge 0}~~\Big | ~~\pi_x^{\ast}L-tE \text{ is a nef line bundle on } \widetilde{X}(x)\right\}.$$
	\end{lemma}

	This definition enables us to conclude that $\varepsilon(X,L,x)\leq \sqrt{L^2}$. If $\varepsilon(X,L,x)<\sqrt{L^2}$ then there exists a reduced and irreducible curve $C$ passing through $x$ such that $\varepsilon(X,L,x)=\frac{L\cdot  C}{\text{mult}_xC}$. Such curves are called \textit{Seshadri curves} of $L$ at $x$. If $\widetilde{C}$ is the strict transform of a Seshadri curve $C$ on $\widetilde{X}(x)$, then by the Hodge index theorem we have $\widetilde{C}^2<0$. 
	
	A reduced irreducible curve $C$ on $X$ is said to be a \textit{$(-m)$-curve} if $C^2=-m$ and $K_X\cdot  C=m-2$, where $K_X$ is the canonical line bundle of $X$.
	
	We now introduce the geometric setting in which our study will take place.
	Let $\pi:X_r \rightarrow\proj^2$ be the blow up of $\proj^2$ at $r$ distinct points $p_1,p_2,\dots,p_r$. Define $H$ to be the pullback of a line not passing through any $p_i$  and let $E_i=\pi^{-1}(\{p_i\})$ denote the exceptional divisor corresponding to $p_i$, for all $1\leq i \leq r$ . Then the linear equivalence classes of $H,E_1,E_2,\dots,E_r$ form a basis for the Picard group $\text{Pic}(X_r)$. The following define intersection product on Pic$(X_r)$: $$H^2=-E_i^2=1 \text{ and } H\cdot  E_i=E_i\cdot  E_j=0 \text{ for all }i\neq j.$$ With this basis, the canonical bundle of $X_r$ is $K_{X_r}:=-3H+E_1+E_2+\dots+E_r$. As the Seshadri constant at a point $x$ on $X_r$ has an alternative definition as defined above, we will often consider the blow up $\pi_x: \Xr \to X_r$ at $x$. By abuse of notation, $E_i$ together with $E_x:= \pi_x^{-1}(\{x\})$ will denote the exceptional divisors for $\Xr$.

	To study the Seshadri constant it is essential to look at linear systems on $\proj^2$ as well as on  $X_r$.	Consider the complete linear system $\mathcal{L}=|dH-m_1E_1-m_2E_2-\dots-m_rE_r|$ for some non-negative integers $m_1,m_2,\dots,m_r$, and $d\geq 1$ on a blow up  $X_r$ of $\proj ^2$ at very general points $p_1,\dots,p_r$. The \textit{expected dimension} of $\mathcal{L}$ is denoted by edim $\mathcal{L}$ and defined as $$  \text{edim } \mathcal{L}=\max \Bigg \{ {d+2 \choose 2} -\displaystyle \sum_{i=1}^{r} {m_i+1 \choose 2}-1, -1\Bigg \}.$$
	It is easy to see that dim $\mathcal{L}\ge$ edim $\mathcal{L}$. We say that $\mathcal{L}$ is \textit{non-special} if dim $\mathcal{L}=\text{edim }\mathcal{L}$ and \textit{special} otherwise.
	
	Segre-Harbourne-Gimigliano-Hirschowitz (SHGH) conjecture is one of the main conjectures in this direction.  Even though there are several versions of this conjecture proposed by Segre\cite{Ser}, Harbourne\cite{Har1}, Gimigliano\cite{Gim} and Hirschowitz\cite{Hir}, we only state the  version proposed by Segre.
	\begin{conjecture}[SHGH Conjecture]\label{SHGH}
		If $\mathcal{L}$ is special, then every divisor in $\mathcal{L}$ is non-reduced.
	\end{conjecture}
	It is known that the SHGH conjecture is true when $r\leq 8$. We recall the following conjecture known as  the Strong SHGH conjecture or Generalized SHGH conjecture \cite[Conjecture 3.6]{Far}, which is stronger than the SHGH conjecture. This was proposed by \L. Farnik, K. Hanumanthu, J. Huizenga, D. Schmitz and T. Szemberg in 2020 in order to exhibit the existence of ample line bundles with an irrational Seshadri constant at every point of a general blow up of $\proj^2$.
	
	\begin{conjecture}[Generalized SHGH or Strong SHGH]\label{SSHGH}
		Let $X_r$ be a blow up of $r$ very general points of $\proj^2$ and let $d\geq 1$, $t\geq 1$ and $m_1,m_2,\dots,m_r \geq 0$ be integers such that
		$$   {d+2 \choose 2} -\displaystyle \sum_{i=1}^{r} {m_i+1 \choose 2} \leq \max \Bigg \{ {t+1 \choose 2}-2,0 \Bigg \}.$$
		Then any curve $C \in |dH-m_1E_1-\dots-m_rE_r|$ which has a point of multiplicity t is non-reduced.
	\end{conjecture}
	In particular, the case $t=1$ is equivalent to the SHGH conjecture.
	The following conjecture follows from the SHGH conjecture.
	\begin{conjecture}[Weak SHGH]\label{WSHGH}
		Let $X_r$ be a blow up of $r$ very general points of $\proj^2$. If $C$ is an irreducible and reduced curve on $X_r$ such that $C^2<0$, then $C$ is a $(-1)$-curve.
	\end{conjecture}
	For $r \geq 10$, Conjecture \ref{SHGH} and Conjecture \ref{WSHGH} are open and Conjecture \ref{SSHGH} is open even for $r\geq 1$.

	\section{Blow ups of points on a curve}
	Let $X_r$ denote a blow up of $\proj^2$ at $r$ points on a curve of degree $d$ on $\proj ^2$. In this section, we will exhibit a lower bound for the Seshadri constant of ample line bundles at any point of $X_r$.

	\begin{theorem}\label{dcurve}
		Let $p_1,p_2,\dots,p_r$ be r distinct points lying on a curve $C\subseteq \mathbb{P}^2$ of degree $d\leq 3$. Let $X_r$ be the blow up of $\mathbb{P}^2$ at $p_1,p_2,\dots ,p_r$. Then for ample line bundle $L$ and for all points $x\in X_r$ we have,
		$$
		\varepsilon(X_r,L,x) \geq 
		\begin{cases}
			1 \text{ if $d=1,2$ or $d=3$ and $L\cdot (-K_{X_r})\geq 2$}\\
			\frac{1}{2} \text{ if $d=3$ and $L\cdot (-K_{X_r})=1.$}
		\end{cases}
		$$
		\begin{proof}
			By \cite[Theorem III.1 (a)]{Har2}, an ample line bundle on a smooth anti-canonical rational surface $X$ is base point free if $L\cdot (-K_X)\geq 2$. It is well-known that for an ample base point free line bundle $L$ on a surface $X$, $\varepsilon(X,L,x)\geq 1$ for every $x\in X$ \cite[Example 5.1.18]{Laz}. First, we consider the case $d=3$. If $L\cdot (-K_{X_r})\geq 2$ then $\varepsilon(X_r,L,x)\geq 1$. If $L\cdot (-K_{X_r})=1$,  $2L$ will be an ample base point free line bundle, which gives $\varepsilon(X_r,2L,x)\geq 1$. This is equivalent to $\varepsilon(X_r,L,x)\geq \frac{1}{2}$.
			
		Now assume $d=1$ or $d=2$. Let $L$ be an ample line bundle on $X_r$ and let $\widetilde{C}$ be the strict transform of $C$ on $X_r$. Since $L$ is ample we have $L\cdot \widetilde{C}\geq 1$ and $L\cdot H\geq 1$. As $-K_{X_r}=\widetilde{C}+(3-d)H$ and $d\leq 2$, it follows that $L\cdot (-K_{X_r}) \geq 2$.
		\end{proof}
	\end{theorem}
	Now we consider the case of points lying on a curve of arbitrary degree.
	\begin{proposition}
		Let $C$ be a smooth curve of degree $d$ on $\proj^2$ and $X_r$ the blow up of $r$ distinct smooth points $p_1,p_2,\dots,p_r \in C$. If $L=eH-n(E_1+\dots + E_r)$ is an ample line bundle on $X_r$ then, $$\varepsilon(X_r,L)\geq \frac{1}{r}, ~\; ~~\; \forall r \geq d^2+1.$$
		\begin{proof}
			Let \ $\widetilde{C}=dH-E_1-\dots-E_r$ be the strict transform of $C$. As above, it suffices to show that $rL$ is globally generated. By \cite[Theorem 3.2]{Han} for $r\geq d^2+1$, $rL$ is globally generated if $(re+3)d>r(rn+1)$. This is equivalent to $r(ed-rn)>r-3d$. This is true since $L\cdot \widetilde{C} =ed-rn \geq 1$.
		\end{proof}
	\end{proposition}

	\section{Blow ups of $r\leq 8$ distinct points}
	In this section, we consider blow ups of $\proj^2$ at $r\leq 8$ points and give lower bounds for Seshadri constants. These partially follow from Theorem \ref{dcurve}. Throughout this section, we will assume that $r\leq 8$ and $X_r$ denotes the blow up of $\mathbb{P}^2$ at $r$ distinct points. We recall the following proposition.
	\begin{proposition}\cite[Proposition 2.6]{KCSA} \label{SC8}
		Let $X_r$ be the blow up of $r\ge 1$ distinct points on $\proj^2$ and let $L$ denote an ample line bundle on $X_r$. Assume $r\leq 8$, and for $x\in X_r$, let $\Xr$ denote the blow up of $X_r$ at $x$. Then
		$$\varepsilon(X_r,L,x)=\inf \bigg \{ \frac{L\cdot C}{\rm{mult}_x C}~~ \bigg | ~~\widetilde{C}\in \Gamma \bigg \},$$
		where $C$ is a curve on $X_r$ and $\Gamma $ is the set of all $(-1)$-curves, $(-2)$-curves, fixed components of $\big |-K_{\Xr}\big |$ and reduced irreducible curves in the linear system $\big |-K_{\Xr}\big |$ on $\Xr$.
	\end{proposition}
	\begin{remark}
		If $C$ is a Seshadri curve on $X_r$ for an ample  line bundle at a point $x\in X_r$, then $\widetilde{C}^2<0$, where $\widetilde{C}$ is the strict transform of $C$ on $\Xr$. Note that by \cite[Proposition 4.1 (ii)]{Sano}, for $r\leq 8$, $\widetilde{C}$ is either a $(-1)$-curve or a $(-2)$-curve or a fixed component of $\big |-K_{\Xr}\big |$.
	\end{remark}	
	\begin{proposition}\label{r7geq1}
		Let $X$ be a smooth anti-canonical surface (i.e., $-K_X$ is an effective divisor). If $(-K_X)^2\geq 2$ we have,
		$$\varepsilon(X)\geq 1.$$ 
		\begin{proof}
			Let $L$ be an ample line bundle on $X_r$. Then as an application of the Hodge index theorem we have,
			$$(-K_X\cdot L)^2\geq (-K_X)^2L^2\geq 2.$$
			So $-K_X\cdot L \geq \sqrt{2}$. Also note that $-K_X\cdot L$ is a positive integer. Thus we have $-K_X\cdot L\geq 2$. Therefore by \cite[Theorem III.1 (a)]{Har2}, $L$ is base point free. So $\varepsilon(X_r,L,x)\geq 1$ for all $x \in X_r$. 
		\end{proof}
	\end{proposition}

	\begin{corollary}
		For $r\leq 7$ we have,
		$$\varepsilon(X_r)\geq 1.$$
	\end{corollary}
	
	Next, we consider the $r=8$ case. First we prove the following lemma.
	\begin{lemma}\label{-K}
		Let $X_8$ be a blow up of $\proj^2$ at 8 distinct points. If $L$ is an ample line bundle with $-K_{X_8}\cdot L =1$ then $L=-K_{X_8}$.
		\begin{proof}
			Let $L=eH-n_1E_1-\dots-n_8E_8$ be an ample line bundle with $-K_{X_8}\cdot L =1$, i.e., $3e-\displaystyle \sum_{i=1}^{8}n_i=1$. If $L^2\geq 2$, the Hodge index theorem gives $$(L\cdot -K_{X_8})^2\geq L^2(-K_{X_8})^2\geq 2,$$
			which contradicts the hypothesis. So $L^2=1$, i.e., $e^2-\displaystyle \sum_{i=1}^8 n_i^2=1$.
			So,
			$$9e^2-6e+1=(3e-1)^2=\bigg (\sum_{i=1}^8 n_i\bigg )^2\leq 8\sum_{i=1}^8 n_i^2=8e^2-8.$$
			This reduces to $e^2-6e+9\leq 0$ which happens only if $e=3$. Since $L$ is ample we have $n_i>0$ and $3=e>n_i+n_j$ for all $i\neq j$. This gives $n_i=1$ for all $1\leq i\leq 8$. So $L=-K_{X_8}$.
		\end{proof}
	\end{lemma}
	\begin{theorem}\label{X8}
		Let $X_8$ be a blow up of $\proj^2$ at 8 distinct points $p_1,p_2,\dots ,p_8$.
		Let $L$ be an ample line bundle on $X_8$. Suppose that $\varepsilon(X_8,L,x)<1$ for some $x\in X_8$. Then the following  hold:
		\begin{enumerate}
			\item $p_1,p_2,\dots ,p_8$ are smooth points on an irreducible singular cubic,
			\item $L=-K_{X_8}$, and
			\item $\varepsilon(X_8,L,x)=\frac{1}{2}$.
		\end{enumerate} 
		\begin{proof}
			Suppose that $ \varepsilon(X_8,L,x) < 1$ for an ample line bundle $L$ at some $x\in X_8$. By Proposition \ref{SC8}, the Seshadri constant is achieved by a curve $C$ on $X_8$ whose strict transform on $\widetilde{X}_8(x)$ is  a $(-1)$-curve or a $(-2)$-curve or a fixed component of  $|-K_{\widetilde{X}_8(x)}|$ or a curve in the linear system $|-K_{\widetilde{X}_8(x)}|$. Let $C=dH-m_1E_1-\dots-m_8E_8$ and $t=\text{mult}_xC >1$, i.e.,
			$$\varepsilon(X_8,L,x)=\frac{L\cdot C}{\text{mult}_xC}=\frac{L\cdot C}{t}.$$
			First, we show that the strict transform of $C$ on $\widetilde{X}_8(x)$ is \textit{not} a $(-1)$-curve or a $(-2)$-curve. If $C$ is such a curve, we have $C^2=t^2-i$, where $i=1,2$. Since
			$$\varepsilon(X_8,L,x)=\frac{L\cdot C}{\text{mult}_xC}<1,$$
			we have $L\cdot C \leq t-1$. By the Hodge index theorem, $(L\cdot C)^2\geq L^2C^2 \geq t^2-i$. Putting these together, we get
			$$(t-1)^2\geq (L\cdot C)^2\geq t^2-i.$$
			This is equivalent to $t\leq \frac{1+i}{2}$, where $i=1$ or $i=2$. This implies that $t=1$, which is a contradiction.
			
			So $C$ is a singular curve whose strict transform is either a fixed component of $|-K_{\widetilde{X}_8(x)}|$ or a curve in the linear system $|-K_{\widetilde{X}_8(x)}|$. Since $C$ is a singular curve, we must have $d=3$. Hence $\widetilde{C}\in |-K_{\widetilde{X}_8(x)}|$. So $p_1,\dots ,p_8$ belong to the image of $C$ on $\proj^2$. As $C$ is irreducible, it is clear that $\text{mult}_xC =2$. Since $\varepsilon(X_8,L,x)<1$, the only possibility is 
			$$\varepsilon(X_8,L,x)=\frac{L\cdot C}{2}=\frac{1}{2}.$$
			Hence $L\cdot -K_{X_8}=1$ and by Lemma \ref{-K} this happens only if $L=-K_{X_8}$.
		\end{proof}
	\end{theorem}
	
	\begin{remark}
		A part of Theorem \ref{X8} has an alternate proof using Theorem \ref{dcurve}. Note that there is a cubic passing through any distinct 8 points in $\proj^2$. Theorem \ref{dcurve} gives that $\frac{1}{2}$ is a lower bound for the Seshadri constant of $X_8$. More precisely, using Lemma \ref{-K} we can see that $\varepsilon(X_8,L)\ge 1$ for all ample line bundles $L\ne -K_{X_8}$ and $\varepsilon(X_8,-K_{X_8})\ge \frac{1}{2}$ provided $-K_{X_8}$ is  an ample line bundle.
	\end{remark}
	
	\begin{example}\label{example1}
		By Theorem \ref{X8}, the least Seshadri constant that can be achieved on $X_8$ is $\frac{1}{2}$. We now give an example in which $\frac{1}{2}$ is achieved at some point on a del Pezzo surface. Consider the smooth rational surface obtained by blowing up 8 smooth points $p_1,p_2,\dots p_8$ on a singular irreducible cubic on $\proj^2$ such that no three lie on a line, no six lie on a conic, and if a singular cubic passes through all the points then $p_i$ are smooth points of that cubic for $1\le i \le 8$. Hence $-K_{X_8}$ is ample by \cite[Th\'eor\`eme 1]{Dem}. Let $x$ denote the point in $X_8$ corresponding to the singular point on the cubic. Let $C$ denote the strict transform of the singular cubic containing the 8 points. From Theorem \ref{X8} it is clear that $\varepsilon(X_8,-K_{X_8},x)\geq \frac{1}{2}$. Note that $\text{mult}_xC=2$ and $-K_{X_8}\cdot C=1$. These give us:
		$$\varepsilon(X_8,-K_{X_8},x)=\frac{1}{2}.$$
	\end{example}

	\begin{proposition}
		For $r\leq 5$, let $X_r$ denote the blow up of $r$ distinct points in $\proj^2$, then $\varepsilon(X_r,L,x)$ is an integer for all ample line bundles $L$ on $X_r$ and for all $x\in X_r$.
		\begin{proof}
			It is known that the number of $(-1)$ and $(-2)$-curves on $\Xr$ is finite for $r\leq 7$. Also for $r\leq 5$, we can see that curves whose strict transforms are $(-1)$ or $(-2)$-curves on $\Xr$ are smooth at $x$. The curves on $X_r$ whose strict transform is a fixed component of $|-K_{\Xr}|$ and is not an element of $|-K_{\Xr}|$ are also smooth. It is also clear that a curve $C$ on $X_r$ such that $\widetilde{C}$ is an element of $|-K_{\Xr}|$ has non-negative self intersection if $r\leq 5$. Note the fact that if $C$ is a Seshadri curve, its strict transform $\widetilde{C}$ has negative self intersection on $\Xr$. This shows that all possible Seshadri curves are smooth at $x$ when $r\leq 5$. Hence the Seshadri constant of every ample line bundles at every point of $X_r$ is an integer when $r\leq 5$.
		\end{proof}
	\end{proposition}

	\begin{remark}
		It is shown in \cite[Th\'eor\`eme 1.3]{Bro} that there is a del Pezzo surface $X_6$ which is the blow up of six points on $\mathbb{P}^2$ such that no three lie on a line, no six lie on a conic, with fractional Seshadri constant. More precisely, we have $\varepsilon(X_6,-K_{X_6},x)=\frac{3}{2}$ for some $x\in X_6$.
	\end{remark}
	
	\section{General blow ups of $\proj ^2$}\label{GenP2}
	We proceed with a similar setting as in the above section. Furthermore throughout this section we assume that $ p_1,p_2,\dots,p_r$ are very general points in $\proj ^2$, where $r\ge 1$. In this section our goal is to find some lower bounds on Seshadri constants on $X_r$ assuming the Strong SHGH Conjecture (Conjecture \ref{SSHGH}) is true.
	
	Following result classify all the possibilities of a self intersection negative curve on $\Xr$ for any $x\in X_r$.
	\begin{lemma}\label{neg}
		For any point $x\in X_r$, let $\Xr$ be the blow up of $X_r$ at x. Assume that Conjecture \ref{SSHGH} is true. Let $C$ be a reduced, irreducible curve on $X_r$ and let $\widetilde{C}$ denote its strict transform on $\Xr$. If $\widetilde{C}^2 < 0$ then $\widetilde{C}$ is one of the following holds:
		\begin{enumerate}
			\item $\widetilde{C}^2 = -1$ and $K_{\Xr}\cdot  \widetilde{C}= -1$,
			\item $\widetilde{C}^2 = -1$ and $K_{\Xr}\cdot  \widetilde{C}= 1$, 
			\item $\widetilde{C}^2 = -1$ and $K_{\Xr}\cdot  \widetilde{C}= 3$,
			\item $\widetilde{C}^2 = -2$ and $K_{\Xr}\cdot  \widetilde{C}= 0$,
			\item $\widetilde{C}^2 = -2$ and $K_{\Xr}\cdot  \widetilde{C}= 2$,
			\item $\widetilde{C}^2 = -3$ and $K_{\Xr}\cdot  \widetilde{C}= 1$.
		\end{enumerate}
		\begin{proof}
			Let $C=dH-m_1E_1-\dots-m_rE_r$ and $t=\text{mult}_x C$. So $\widetilde{C} = dH-m_1E_1-\dots-m_rE_r-tE_x$, where $E_x$ is the exceptional divisor on $\Xr$. If $d=0$, $\widetilde{C}=E_i$ or $\widetilde{C}=E_i-E_x$. Both are of the required type.
			
			 Now let $d\ge 1$. Since $C$ is reduced, Conjecture \ref{SSHGH} implies
			$$   {d+2 \choose 2} -\displaystyle \sum_{i=1}^{r} {m_i+1 \choose 2} > {t+1 \choose 2}-2.$$ After some computations and rearrangements, the above inequality reduces to
			$$\frac{d^2- \displaystyle \sum_{i=1}^{r} m_i^2 - t^2 +3d -\sum_{i=1}^{r} m_i-t}{2}+1>-2,$$
			which is same as
			\begin{equation}\label{RRTgrt2}
			\frac{\widetilde{C}^2-K_{\Xr}\cdot\widetilde{C}}{2}+1>-2.
			\end{equation}
			The Riemann-Roch theorem gives,
			$$h^0(\Xr,\widetilde{C})-h^1(\Xr,\widetilde{C})+h^2(\Xr,\widetilde{C})=\frac{\widetilde{C}^2-K_{\Xr}\cdot  \widetilde{C}}{2}+1.$$
			Since $\widetilde{C}$ is an effective divisor on $\Xr$, we get $h^2(\Xr,\widetilde{C})=0$. 
			Also $\widetilde{C}^2<0$ implies $h^0(\Xr,\widetilde{C})= 1$. Then \eqref{RRTgrt2} reduces to:
			\begin{equation}\label{eeqn}
				 -h^1(\Xr,\widetilde{C})=\frac{\widetilde{C}^2-K_{\Xr}\cdot  \widetilde{C}}{2}>-3.
			\end{equation}
			
			Hence $h^1(\Xr,\widetilde{C})\leq 2.$
			We consider the following three cases:
			\begin{itemize}
				\item[Case 1: ]
				$h^1(\Xr,\widetilde{C}) =0$. \\
				From \eqref{eeqn} we get, $ h^1(\Xr,\widetilde{C}) =\frac{K_{\Xr}\cdot  \widetilde{C}-\widetilde{C}^2}{2}=0 $. i.e., $K_{\Xr}\cdot \widetilde{C}=\widetilde{C}^2$.\\
				Now from adjunction formula we see that, $$K_{\Xr}\cdot  \widetilde{C}+\widetilde{C}^2 \geq -2.$$
				Since $\widetilde{C}^2 <0$, we get $K_{\Xr}\cdot  \widetilde{C}=\widetilde{C}^2=-1.$
				
				\item[Case 2: ]
				$h^1(\Xr,\widetilde{C}) =1$. \\
				Then $K_{\Xr}\cdot  \widetilde{C}-\widetilde{C}^2=2$. As above, using the adjunction formula and $\widetilde{C}^2<0$, we conclude that $\widetilde{C}^2 = -1$ and $K_{\Xr}\cdot  \widetilde{C}= 1$ or $\widetilde{C}^2 = -2$ and $K_{\Xr}\cdot  \widetilde{C}= 0$.
				\item[Case 3: ] $h^1(\Xr,\widetilde{C}) =2$. \\
				Similarly as in above cases, we conclude that types (3) or (5) or (6) can occur.
			\end{itemize}
		\end{proof}
	\end{lemma}
	In fact Lemma \ref{neg} says that $C^2\ge -3$ for all reduced irreducible curve $C$ on $\Xr$. The following is the main result of this section.
	\begin{theorem}\label{scgeq1}
		Assume that Conjecture \ref{SSHGH} is true. Let $\pi : X_r \to \proj ^2$ be the blow up of $\proj ^2$ at very general points $p_1,p_2,\dots ,p_r$. If $\varepsilon(X_r,L,x)<1$ for some ample line bundle $L$ and for some $x\in X_r$, then $\varepsilon(X_r,L,x)=\frac{1}{2}$. In particular, $$\varepsilon(X_r)\ge \frac{1}{2}.$$
		\begin{proof}
			Let $L$ be an ample line bundle on $X_r$ with $\varepsilon(X_r,L,x) < 1$ for some $x\in X_r$. Then $\varepsilon(X_r,L,x)$ is achieved by an irreducible and reduced curve $C\subseteq X_r$, i.e.,
			$$\varepsilon(X_r,L,x)=\frac{L\cdot C}{\text{mult}_xC}.$$ 
						
			We know that $\widetilde{C}^2<0$, where $\widetilde{C}$  is the strict transform of $C$. Let $t=\text{mult}_xC$. Then $\widetilde{C}^2=C^2-t^2$. Since $\varepsilon(X_r,L,x)<1$ we have $t>1$ and $L \cdot C \leq t-1$. By Lemma \ref{neg} we have, $0> \widetilde{C}^2 \geq -3$. Hence $C^2=t^2-i$, where $i=1 $ or $2$ or $3$. 
			
			Now using the Hodge index theorem,
			$$(t-1)^2\geq (L \cdot C)^2 \geq L^2C^2 \geq t^2-i.$$
			So $(t-1)^2\geq  t^2-i$, which reduces to $t\leq \frac{1+i}{2}$ for $i=1,2,3$. Hence $t=1$ or $t=2$. But $\varepsilon(X_r,L,x)<1$ forces $t=2$ and hence $\varepsilon(X_r,L,x)=\frac{1}{2}$.
		\end{proof}
	\end{theorem}
	\begin{remark}
		From the above theorem, it follows that if Conjecture \ref{SSHGH} is true then $\varepsilon(X_r,x)$ and $\varepsilon(X_r)$ cannot be zero for blows ups of a finite number of very general points in $\proj ^2$. A weaker assumption than Conjecture \ref{SSHGH} gives the same conclusion, as we show below.
	\end{remark}
	\begin{conjecture}\label{NCon}
		Let $X_r$  be as in Theorem \ref{scgeq1} and let $d\geq 1$, $t\geq 2$ and $m_1,m_2,\dots,m_r \geq 0$ be integers such that 
		$$   {d+2 \choose 2} -\displaystyle \sum_{i=1}^{r} {m_i+1 \choose 2} \leq {t \choose 2}.$$
		Then any curve $C \in |dH-m_1E_1-\dots-m_rE_r|$ which has a point of multiplicity t is non-reduced.
	\end{conjecture}
	Note that for $t\ge 2$, ${t+1 \choose 2}-2 \geq {t \choose 2}$. Hence Conjecture \ref{NCon} follows from Conjecture \ref{SSHGH}.
	\begin{proposition}
		Let $X_r$ be as in Theorem \ref{scgeq1}. If Conjecture \ref{NCon} is true then $\varepsilon(X_r)\geq \frac{1}{2}$.
		\begin{proof}
			Let $\pi_x,\Xr,E_x$ be as above. Let $L=eH-n_1E_1-n_2E_2-\dots -n_rE_r$ be an ample line bundle on $X_r$ and $x\in X_r$ be a point.
			To show that $\varepsilon(X_r,L,x)\geq \frac{1}{2}$, it is enough to show that
			$\pi_x^{\ast}L-\frac{1}{2}E_x$ is nef, i.e., 	$2\pi_x^{\ast}L-E_x$ is nef. Equivalently, for any reduced, irreducible curve $\widetilde{C}=dH-m_1E_1-\dots-m_rE_r-tE_x$ on $\Xr$, we need to show that
			$$2(de-m_1n_1-\dots-m_rn_r)-t\geq 0.$$ 
			
			Since $L$ is ample, we have $L\cdot C=de-m_1n_1-\dots-m_rn_r\ge 1$. Hence if $t=1$ clearly $2(de-m_1n_1-\dots-m_rn_r)-t\geq 0.$ Also if $d=0$ we have, $\widetilde{C}=E_i \text{ or } \widetilde{C}=E_i-E_x$ clearly $(\pi_x^{\ast}L-\frac{1}{2}E_x)\cdot\widetilde{C}\geq 0$. 
			
			Now assume $d\geq 1$ and $t\ge 2$. Since $\widetilde{C}$ is a reduced curve on $\Xr$, $C=\pi_x(\widetilde{C})$ is a reduced curve on $X_r$. So using Conjecture \ref{NCon} we conclude:
			$$   {d+2 \choose 2} -\displaystyle \sum_{i=1}^{r} {m_i+1 \choose 2} > {t \choose 2}.$$
			
			Consider a very general point $y\in X_r$ and let $p_{r+1}:=\pi(y).$ Therefore, by dimension count, there is a curve $D$ on $\proj^2$ such that $\text{mult}_{p_i}D \geq m_i$ for all $i=1,2,\dots,r$ and $\text{mult}_{p_{r+1}}D \geq t-1$.\\
			By \cite[Theorem]{EinLaz} we have, $\varepsilon(X_r,L,y)\geq 1$. From this we conclude
			$$de-m_1n_1-\dots-m_rn_r-(t-1)\geq 0$$
			and ampleness of $L$ gives
			$$de-m_1n_1-\dots-m_rn_r \geq 1.$$
			Taking the sum of these two inequalities, we get
			$$2(de-m_1n_1-\dots-m_rn_r)-t\geq 0.$$ This gives the required inequality: $$\varepsilon(X_r,L,x)\geq \frac{1}{2}.$$
		\end{proof}
	\end{proposition}
	
	One of the most important questions in the area of Seshadri constants concerns their irrationality, i.e., is there an ample line bundle on a smooth surface so that at some point in the surface, the Seshadri constant of that line bundle is irrational? We propose the following question:
	\begin{question}\label{ques}
		For a given $\frac{p}{q}\in \mathbb{Q}$, is there a positive integer $r$ and a line bundle $L$ on $X_r$ such that $\varepsilon(X_r,L,x)=\frac{p}{q}$ for some $x\in X_r$.
	\end{question}
	Clearly Theorem \ref{scgeq1} implies that if Conjecture \ref{SSHGH} is true, Question \ref{ques} has a negative answer for $\frac{p}{q} <1$ unless $\frac{p}{q}=\frac{1}{2}$. 
	It is enough to look at $\frac{p}{q}$ with $p\geq q$. For this, it is enough to look for line bundles with Seshadri constant $\frac{p}{q}$, where $p$ is a prime number as $\varepsilon(X,nL,x)=n\varepsilon(X,L,x).$
	
	We give below a partial result in this direction. The following result shows all integers are attained as Seshadri constants.
	\begin{proposition}
		For any $r$, there is an ample line bundle $L$ and $x\in X_r$ such that $\varepsilon(X_r,L,x)=1$. Further, every positive integer is achieved as a Seshadri constant for some ample line bundle.
		\begin{proof}
			By \cite[Proposition 3.2]{GGP}, there exists a $d\in \mathbb{N}$ such that $L=dH-E_1-\dots-E_r$ is a very ample line bundle  on $X_r$. Since $L$ is very ample we have $\varepsilon(X_r,L,x)\geq 1$ for all $x\in X_r$.  Since $L\cdot E_i=1$ we have, $\varepsilon(X_r,L,x)\leq 1$ for all  $x \in E_i$ for every $1\leq i \leq r$. So $\varepsilon(X_r,L,x)=1$ for all $x \in \displaystyle \cup_{i=1}^rE_i$.
		\end{proof}
	\end{proposition}

\end{document}